\documentclass[a4paper,12pt,reqno]{article}
\usepackage[T1]{fontenc}
\usepackage[hmargin=1in,vmargin=1in,nohead]{geometry}
\usepackage{textcomp}
\usepackage{amsmath}
\usepackage{graphicx,epsfig,color}
\usepackage{amssymb}
\usepackage{esint}

\makeatletter

\newtheorem{thm}{Theorem}
\newtheorem{lem}[thm]{Lemma}

\newtheorem{prop}[thm]{Proposition}
\newtheorem{rem}{Remark}[section]
\newtheorem{remark}{Remark}[section]

\newtheorem{Def}{Definition}
\DeclareMathOperator{\cg}{\textbf{[}}
\DeclareMathOperator{\cd}{\textbf{]}}
\newcommand{\cro}[1]{\cg {#1} \cd}
\title{ \bf Generalized Fleming-Viot processes with immigration via stochastic flows of partitions}
\author{ \bf Cl\'ement Foucart \\
\\
\emph {Laboratoire de Probabilit\'es et Mod\`eles Al\'eatoires} \\
\emph {Universit\'e Pierre et Marie Curie}\\ 
\emph{4 Place Jussieu- 75252 Paris Cedex 05- France}}
\begin{document}
\maketitle{}
\begin{abstract}
The generalized Fleming-Viot processes were defined in 1999 by Donnelly and Kurtz using a particle model and by Bertoin and Le Gall in 2003 using stochastic flows of bridges. In both methods, the key argument used to characterize these processes is the duality between these processes and exchangeable coalescents. A larger class of coalescent processes, called distinguished coalescents, was set up recently to incorporate an immigration phenomenon in the underlying population. The purpose of this article is to define and characterize a class of probability measure-valued processes called the generalized Fleming-Viot processes with immigration. We consider some stochastic flows of partitions of $\mathbb{Z}_{+}$, in the same spirit as Bertoin and Le Gall's flows, replacing roughly speaking, composition of bridges by coagulation of partitions. Identifying at any time a population with the integers $\mathbb{N}:=\{1,2,...\}$, the formalism of partitions is effective in the past as well as in the future especially when there are several simultaneous births. We show how a stochastic population may be directly embedded in the dual flow. An extra individual $0$ will be viewed as an external generic immigrant ancestor, with a distinguished type, whose progeny represents the immigrants. The "modified" lookdown construction of Donnelly-Kurtz is recovered when neither simultaneous multiple births nor immigration are taken into account. In the last part of the paper we give a sufficient criterion for the initial types extinction. 
\end{abstract}
 \vspace{9pt} \noindent {\bf Key words.}
{Exchangeability}, {duality}, {stochastic flow}, {lookdown construction}, {immigration}, {Fleming-Viot process}, {measure-valued process}, {coalescent}, {partition-valued process}, {coming down from infinity}.
\par \vspace{9pt}
\par \vspace{9pt} \noindent {\bf e-mail.} {clement.foucart@etu.upmc.fr}
\par
\newpage
\section{Introduction}
Originally, Fleming-Viot processes were defined in 1979 in \cite{FV} to model the genetic phenomenon of allele drift. They now form an important sub-class of measure-valued processes, which have received significant attention in the literature. Donnelly and Kurtz, in 1996, established in \cite{DonnKurtz1} a duality between the classical Fleming-Viot process and the Kingman coalescent (defined in \cite{Kingman3}). Shortly after, the class of coalescent processes was considerably generalized by assuming that multiple coagulations may happen simultaneously. We refer to Chapter 3 of \cite{Beres2} for an introduction to the $\Lambda$-coalescents and to the seminal papers \cite{Mohle}, \cite{Schweinsberg} for the definition of the general exchangeable coalescent processes (also called $\Xi$-coalescents). An infinite particle representation set up in 1999, allowed Donnelly and Kurtz, in \cite{DonnKurtz} to define a generalized Fleming-Viot process by duality with the $\Lambda$-coalescent. In 2003, Bertoin and Le Gall started from another point of view and introduced, in \cite{LGB1}, a stochastic flow of bridges which encodes simultaneously an exchangeable coalescent process and a continuous population model, the so-called generalized Fleming-Viot process. Finally, in 2009, Birkner \textit{et al}. in \cite{Birk2} have adapted the same arguments as Donnelly-Kurtz for the case of the $\Xi$-coalescent.
\\

More recently, a larger class of coalescents called distinguished coalescents were defined, see \cite{coaldist}, in order to incorporate immigration in the underlying population. The purpose of this article is to define by duality the class of generalized Fleming-Viot processes with immigration. Imagine an infinite haploid population with immigration, identified by the set $\mathbb{N}$. This means that each individual has at most one parent in the population at the previous generation; indeed, immigration implies that some individuals may have parents outside this population (they are children of immigrants). To encode the arrival of new immigrants, we consider the external integer $0$ as the generic parent of immigrants. We shall then work with partitions of $\mathbb{Z}_{+}$, the so-called \textit{distinguished} partitions. Our approach will draw both on the works of Bertoin-Le Gall and of Donnelly-Kurtz. Namely, in the same vein as Bertoin and Le Gall's article \cite{LGB1}, we define a \textit{stochastic flow of partitions} of $\mathbb{Z}_{+}$, denoted by $(\Pi(s,t), -\infty<s\leq t<\infty)$. The dual flow $(\hat{\Pi}(t), t\geq 0):=(\Pi(-t,0),t\geq 0)$ shall encode an infinite haploid population model with immigration forward in time. Namely, for any individual $i\geq 1$ alive at the initial time $0$, the set $\hat{\Pi}_{i}(t)$ shall represent the descendants of $i$, and $\hat{\Pi}_{0}(t)$ the descendants of the generic immigrant, at time $t$. We stress that the evolution mechanism involves both multiple reproduction and immigration. The genealogy of this population model is precisely a distinguished exchangeable coalescent. Our method is close to the "modified" lookdown construction of Donnelly-Kurtz, but differs from (and simplifies) the generalization given by Birkner et al. in \cite{Birk2}. In the same vein as Fleming and Viot's fundamental article \cite{FV}, we consider the type carried by each individual at the initial time. Denote by $U_{i}$ the type of the individual $i\geq 1$, and distinguish the type of the generic immigrant by fixing $U_{0}=0$. The study of the evolution of frequencies as time passes, leads us to define and study the so-called generalized Fleming-Viot process with immigration (GFVI for short) denoted in the sequel by $(Z_{t}, t\geq 0)$. This process will be explicitly related to the \textit{forward partition-valued process} $(\hat{\Pi}(t),t\geq 0)$. When the process $(\hat{\Pi}(t),t\geq 0)$ is absorbed in the trivial partition $(\{\mathbb{Z}_{+}\}, \emptyset,...)$, all individuals are immigrant children after a certain time. We shall discuss conditions entailing the occurrence of this event. 
\\
\\ 
The rest of the paper is organized as follows. In Section 2 (Preliminaries), we recall some basic facts on distinguished partitions. We introduce a coagulation operator $coag$ and describe how a population may be encoded by distinguished partitions forward in time. Some properties related to exchangeable sequences and partitions are then presented. We recall the definition of an exchangeable distinguished coalescent and define a stochastic flow of partitions using a Poisson random measure on the space of partitions. In Section 3, we study the dual flow and the embedded population. Adding initially a type to each individual, the properties of exchangeability of Section 2 allow us to define the generalized Fleming-Viot process with immigration. We show that any generalized Fleming-Viot process with immigration is a Feller process. Arguing then by duality we determine the generator of any GFVI on a space of functionals which forms a core. In Section 4, we give a sufficient condition for the extinction of the initial types. Thanks to the duality established in Section 3, the extinction corresponds to the coming down from infinity of the distinguished coalescent. 
\section{Preliminaries}
We start by recalling some basic definitions and results about exchangeable distinguished partitions which are developed in \cite{coaldist}. For every $n\geq 0$, we denote by $\cro{n}$ the set $\{0,...,n\}$ and call $\mathcal{P}^{0}_{n}$ the set of partitions of $\cro{n}$. The space $\mathcal{P}^{0}_{\infty}$ is the set of partitions of $\cro{\infty}:=\mathbb{Z}_{+}$. By convention the blocks are listed in the increasing order of their least element and we denote by $\pi_{i}$ the $i$-th block of the partition $\pi$. The first block of $\pi$, which contains $0$, is thus $\pi_{0}$ and is viewed as distinguished. It shall represent descendants of an immigrant. An element of $\mathcal{P}^{0}_{\infty}$ is then called a distinguished partition. We endow the space $\mathcal{P}^{0}_{\infty}$ with a distance, that makes it compact, defined by 
$$d(\pi,\pi')=(1+\max\{n\geq 0; \pi_{|\cro{n}}=\pi'_{|\cro{n}}\})^{-1}.$$
The notation $i\overset{\pi}{\sim}j$ means that $i$ and $j$ are in the same block of $\pi$. An exchangeable distinguished partition is a random element $\pi$ of $\mathcal{P}^{0}_{\infty}$ such that for all permutations $\sigma$ satisfying $\sigma(0)=0$, the partition $\sigma \pi$ defined by \begin{center} $i\overset{\sigma \pi}\sim j$ if and only if $\sigma(i)\overset{\pi}\sim \sigma(j),$ \end{center} has the same law as $\pi$. If $\pi$ is an exchangeable distinguished  partition, the asymptotic frequency  $$|\pi_{i}|:=\underset{n\rightarrow \infty}{\lim}\frac{\#(\pi_{i}\cap \cro{n})}{n}$$ exists for every $i\geq 0$, almost-surely. 
We denote by $|\pi|^{\downarrow}$ the sequence of asymptotic frequencies $(|\pi_{i}|, i\geq 0)$ after a decreasing rearrangement apart from $|\pi_{0}|$. Thanks to Kingman's correspondence, the law of $|\pi|^{\downarrow}$ determines completely that of $\pi$ (see Theorem 2.1 in \cite{coaldist}). 
\\

The operator $coag$ is defined from $\mathcal{P}^{0}_{\infty}\times \mathcal{P}^{0}_{\infty}$ to $\mathcal{P}^{0}_{\infty}$ such that for $(\pi, \pi') \in \mathcal{P}^{0}_{\infty}\times \mathcal{P}^{0}_{\infty}$, the partition $coag(\pi,\pi')$ satisfies for all $i\geq 0$,
$$coag(\pi, \pi')_{i}=\bigcup_{j\in \pi'_{i}}\pi_{j}.$$
The operator $coag$ is Lipschitz-continuous and associative in the sense that for any $\pi, \pi'$ and $\pi''$ in $\mathcal{P}^{0}_{\infty}$
$$coag(\pi, coag(\pi', \pi''))=coag(coag(\pi, \pi'), \pi'').$$  
Moreover, the partition of $\mathbb{Z}_{+}$ into singletons, denoted by $0_{\cro{\infty}}$, may be viewed as a neutral element for the operator $coag$, indeed $$coag(\pi, 0_{\cro{\infty}})=coag(0_{\cro{\infty}},\pi)=\pi.$$
More generally, for every $n\geq 0$, we denote by $0_{\cro{n}}$ the partition of $\cro{n}$ into singletons. It should be highlighted that given two independent exchangeable distinguished partitions $\pi, \pi'$, the partition $coag(\pi,\pi')$ is still exchangeable (see Lemma 4.3 in \cite{coursbertoin} and Remark 2.1 below).
\subsection{Partitions of the population and exchangeability}
We explain how the formalism of partitions may be used to describe a population with immigration as time goes forward. As in the Introduction, imagine an infinite haploid population with immigration, identified by the set $\mathbb{N}$ evolving forward in time. An additional individual $0$ is added and plays the role of a generic immigrant. The model may be described as follows: let $t_{0}<t_{1}$, 
\begin{itemize} 
\item at time $t_{1}$ the families sharing the same ancestor at time $t_{0}$ form a distinguished partition $\pi^{(1)}$. The distinguished block $\pi^{(1)}_{0}$ comprises the children of immigrants,
\item the indices of the blocks of $\pi^{(1)}$ are viewed as the ancestors living at time $t_{0}$. In other words,  for any $j\geq 0$, the block $\pi^{(1)}_{j}$ is the offspring at time $t_{1}$ of the individual $j$ living at time $t_{0}$.
\end{itemize} 
Consider a time $t_{2}>t_{1}$ and denote by $\pi^{(2)}$ the partition of the population at time $t_{2}$ such that the block $\pi^{(2)}_{k}$ comprises the descendants at time $t_{2}$ of the individual $k$, living at time $t_{1}$. Obviously, the set $\bigcup_{k\in \pi_{j}^{(1)}}\pi^{(2)}_{k}$ represents the descendants at time $t_{2}$ of the individual $j$ at time $t_{0}$. Therefore the partition $coag(\pi^{(2)},\pi^{(1)})$ encodes the descendants at time $t_{2}$ of the individuals living at time $t_{0}$.
\\
\\
For any fixed $\pi \in \mathcal{P}^{0}_{\infty}$, define the map $\alpha_{\pi}$ by
\begin{center}
$\alpha_{\pi}(k):=$ the index of the block of $\pi$ containing $k$. 
\end{center}
Thus, in the population above, $\alpha_{\pi^{(1)}}(k)$ corresponds to the ancestor at time $t_{0}$ of the individual $k$ in the population at time $t_{1}$. We have, by definition of the operator $coag$, the key equality
$$\alpha_{coag(\pi,\pi')}=\alpha_{\pi'}\circ \alpha_{\pi}.$$
Therefore, the ancestor living at time $t_{0}$ of the individual $k$ at time $t_{2}$ is $$\alpha_{\pi^{(1)}}\circ \alpha_{\pi^{(2)}}(k)=\alpha_{coag(\pi^{(2)},\pi^{(1)})}(k).$$

We give in the sequel some properties of the map $\alpha_{\pi}$ when $\pi$ is an exchangeable distinguished partition. They will be useful to define a generalized Fleming-Viot process with immigration. We denote by $\mathcal{M}_{1}$ the space of probability measures on $[0,1]$. Let $\rho$ be a random probability measure on $[0,1]$. We say that the exchangeable sequence $(U_{i},i\geq 1)$ has de Finetti measure $\rho$, if conditionally given $\rho$ the variables 
$(U_{i},i\geq 1)$ are i.i.d. with law $\rho$.\\

We make the following key observation.
\begin{lem} \label{exch} Let $(U_{i},i\geq 1)$ be an infinite exchangeable sequence taking values in $[0,1]$, with de Finetti measure $\rho$ and fix $U_{0}=0$. Let $\pi$ be an independent exchangeable distinguished partition. Then the infinite sequence $(U_{\alpha_{\pi}(k)},k\geq 1)$ is exchangeable. Furthermore, its de Finetti measure is $$(1-\sum_{i\geq 0}|\pi_{i}|)\rho+ \sum_{i\geq 1}|\pi_{i}|\delta_{U_{i}}+|\pi_{0}|\delta_{0}.$$ 
\end{lem}
\textit{Proof.} The proof requires rather technical arguments and is given in the Appendix.
\begin{remark}
Let $(U_{i}, i\geq 1)$ be an i.i.d. sequence with a continuous distribution $\rho$. Observing that $$k\overset{\pi}\sim l \Longleftrightarrow U_{\alpha_{\pi}(k)}=U_{\alpha_{\pi}(l)},$$
we can use the previous lemma to recover the (distinguished) paint-box structure of any (distinguished) exchangeable random partition, see Kingman \cite{Kingman1} or Theorem 2.1 in \cite{coursbertoin} for the case with no distinguished block. Moreover, Lemma \ref{exch} yields another simple proof for the exchangeability of $coag(\pi, \pi')$ provided that $\pi$ and $\pi'$ are independent and both exchangeable.
\end{remark} 
\subsection{Distinguished coalescents and flows of partitions}
We start by recalling some basic facts on distinguished coalescents which are developed in Section 3 of \cite{coaldist}. In particular, we recall the definition of a coagulation measure and its decomposition. 
\subsubsection{Distinguished coalescents and coagulation measure} \label{coaldist}
Consider an infinite haploid population model with immigration, the population being identified with $\mathbb{N}$. Recall that a generic immigrant $0$ is added to the population. We denote by $\Pi(s)$ the partition of the current population into families having the same ancestor $s$ generations \textit{earlier}. As explained before, individuals issued from the immigration form the distinguished block  $\Pi_{0}(s)$ of $\Pi(s)$. When some individuals have the same ancestor at a generation $s$, they have the same ancestor at any previous generation. The following statement makes these ideas formal. We stress that time goes backward. An exchangeable distinguished coalescent is a Markov process $(\Pi(t), t\geq 0)$ valued in $\mathcal{P}^{0}_{\infty}$ such that given $\Pi(s)$ $$\Pi(s+t)\overset{Law}{=}coag(\Pi(s), \pi),$$ where $\pi$ is an exchangeable distinguished partition independent of $\Pi(s)$, with a law depending only on $t$. A distinguished coalescent is called standard if $\Pi(0)=0_{\cro{\infty}}$. We emphasize that classical coagulations and coagulations with the distinguished block may happen simultaneously. 
\\
A distinguished coalescent is characterized by a measure $\mu$ on $\mathcal{P}^{0}_{\infty}$, called the coagulation measure, which fulfills the following conditions:
\begin{itemize}
\item $\mu$ is exchangeable, meaning here invariant under the action of the permutations $\sigma$ of $\mathbb{Z}_{+}$ with finite support (i.e permuting only finitely many points), such that $\sigma(0)=0$
\item $\mu(\{0_{\cro{\infty}}\})=0$ and for all $n\geq 0$, $\mu(\pi \in \mathcal{P}^{0}_{\infty}: \pi_{|\cro{n}}\neq 0_{\cro{n}})<\infty$.
\end{itemize}
To be more precise, let $(\Pi(t), t\geq 0)$ be a distinguished coalescent; the coagulation measure $\mu$ is defined from the jump rates of the restricted processes $(\Pi_{|\cro{n}}(t), t\geq 0)$ for $n\geq 0$. For every $\pi \in \mathcal{P}^{0}_{n}$, let $q_{\pi}$ be the jump rate of $(\Pi_{|\cro{n}}(t), t\geq 0)$ from $0_{\cro{n}}$ to $\pi$ and $\mathcal{P}^{0}_{\infty, \pi}$ be the set $$\mathcal{P}^{0}_{\infty, \pi}:=\{\pi'\in \mathcal{P}^{0}_{\infty}; \pi'_{|\cro{n}}=\pi\}.$$ We have by definition $\mu(\mathcal{P}^{0}_{\infty, \pi})=q_{\pi}$. We will denote by $\mathcal{L}_{n}^{*}$ the generator of the continuous Markov chain $(\Pi_{|\cro{n}}(t), t\geq 0)$. Let $\phi$ be a map from $\mathcal{P}^{0}_{n}$ to $\mathbb{R}$ and $\pi \in \mathcal{P}^{0}_{n}$, then
$$\mathcal{L}_{n}^{*}\phi(\pi)=\sum_{\pi'\in \mathcal{P}^{0}_{n}}q_{\pi'}[\phi(coag(\pi,\pi'))-\phi(\pi)].$$
Conversely for any coagulation measure $\mu$, by the same arguments as in \cite{coursbertoin} for the genuine coalescents, a distinguished coalescent with coagulation measure $\mu$, is constructed using a Poisson random measure on the space $\mathbb{R}_{+}\times\mathcal{P}^{0}_{\infty}$ with intensity $dt\otimes\mu$ (see Proposition 3.2 in \cite{coaldist}).  We mention that Theorem 3.1 in \cite{coaldist} yields a decomposition of a coagulation measure into a "Kingman part" and a "multiple collisions part". Let $\mu$ be a coagulation measure, then there exist $c_{0}, c_{1}$ non-negative real numbers and a measure $\nu$ on 
$$\mathcal{P}_{\textbf{m}}:=\left\{s=(s_{0},s_{1},...);\ s_{0}\geq 0, s_{1}\geq s_{2}\geq ...\geq 0, \sum_{i\geq 0}s_{i}\leq 1\right \}$$ 
such that
$$\mu=c_{0}\sum_{1\leq i} \delta_{K(0,i)}+c_{1}\sum_{1\leq i<j} \delta_{K(i,j)}+\int_{s\in \mathcal{P}_{\textbf{m}}}\rho_{s}(.)\nu(ds)$$
where $K(i,j)$ is the simple partition (meaning with at most one non-singleton block) with doubleton $\{i,j\}$ and $\rho_{s}$ denotes the law of an $s$-distinguished paint-box (see Definition 2.3 in \cite{coaldist}). The measure $\nu$ satisfies the condition $$\int_{\mathcal{P}_{\textbf{m}}}(s_{0}+\sum_{i\geq 1}s_{i}^{2})\nu(ds)<\infty.$$ 
\subsubsection{Stochastic flow of partitions}
We define a stochastic flow of partitions and give a construction from a Poisson random measure. The following definition may be compared with that of Bertoin and Le Gall's flows \cite{LGB1}. 
\begin{Def} A flow of distinguished partitions is a collection of random variables $(\Pi(s,t), -\infty<s\leq t<\infty)$ valued in $\mathcal{P}^{0}_{\infty}$ such that:
\begin{itemize}
\item[(i)] For every $t\leq t'$, the distinguished partition $\Pi(t,t')$ is exchangeable with a law depending only on $t'-t$. 
\item[(ii)] For every $t<t'<t''$, $\Pi(t,t'')=coag(\Pi(t,t'), \Pi(t',t''))$ almost surely
\item[(iii)] if $t'_{1}<t'_{2}<...<t'_{n}$, the distinguished partitions $\Pi(t'_{1},t'_{2}),...,\Pi(t'_{n-1},t'_{n})$ are independent.
\item[(v)] $\Pi(0,0)=0_{\cro{\infty}}$ and $\Pi(t,t')\rightarrow 0_{\cro{\infty}}$ in probability when $t'-t \rightarrow 0$.
\end{itemize}
\end{Def}
The process $(\Pi(t), t\geq 0):=(\Pi(0,t), t\geq 0)$ is by definition a distinguished exchangeable coalescent. Given a coagulation measure $\mu$, we introduce and study next a stochastic flow of partitions constructed from a Poisson random measure on $\mathbb{R}\times \mathcal{P}^{0}_{\infty}$ with intensity $dt\otimes \mu$. Instead of composing bridges as Bertoin and Le Gall in \cite{LGB1}, we coagulate directly the partitions replacing thus the operator of composition by the operator $coag$. For any partitions $\pi^{(1)},...,\pi^{(k)}$, we define recursively the partition $coag^{k}(\pi^{(1)},...,\pi^{(k)})$ by $coag^{0}=0_{\cro{\infty}}$, $coag^{1}(\pi^{(1)})=\pi^{(1)}$ and for all $k\geq 2$, 
\begin{align*}
coag^{k}(\pi^{(1)},...,\pi^{(k)})&:=coag \left(coag^{k-1}\left(\pi^{(1)},...,\pi^{(k-1)}\right),\pi^{(k)}\right)\\
&=coag^{k-1}\left(\pi^{(1)},...,\pi^{(k-2)},coag\left(\pi^{(k-1)},\pi^{(k)}\right)\right).
\end{align*}
Introduce a Poisson random measure $\mathcal{N}$ on $\mathbb{R}\times \mathcal{P}^{0}_{\infty}$ with intensity $dt\otimes \mu$ and for each $n\in \mathbb{N}$, let $\mathcal{N}_{n}$ be the image of $\mathcal{N}$ by the map $\pi \mapsto \pi_{|\cro{n}}$. The condition $\mu(\pi_{|\cro{n}}\neq 0_{\cro{n}})<\infty$ ensures that for all $s<t$ there are finitely many atoms of $\mathcal{N}_{n}$ in $]s,t]\times \mathcal{P}^{0}_{n}\setminus{\{0_{\cro{n}}\}}$. We denote by $\{(t_{1}, \pi^{(1)}),(t_{2}, \pi^{(2)}),..., (t_{K}, \pi^{(K)})\}$ these atoms with $K:=\mathcal{N}_{n}(]s,t]\times \mathcal{P}^{0}_{n}\setminus \{0_{\cro{n}}\})$ and define
$$\Pi^{n}(s,t):=coag^{K}\left(\pi^{(1)},...,\pi^{(K)}\right).$$
It remains to establish the compatibility of the sequence of random partitions $(\Pi^{n}(s,t), n\in \mathbb{N})$, which means that for all $m\leq n$, $\Pi^{n}_{|\cro{m}}(s,t)=\Pi^{m}(s,t)$. All non-trivial atoms of $\mathcal{N}_{m}$ are plainly non-trivial atoms of $\mathcal{N}_{n}$, and moreover the compatibility property of the operator $coag$ with restrictions implies that
$$coag^{K}\left(\pi^{(1)},...,\pi^{(K)}\right)_{|\cro{m}}=coag^{K}\left(\pi^{(1)}_{|\cro{m}},...,\pi^{(K)}_{|\cro{m}}\right).$$
Two cases may occur, either $\pi^{(i)}_{|\cro{m}}=0_{\cro{m}}$ and does not affect the coagulation, or $\pi^{(i)}_{|\cro{m}}\neq 0_{\cro{m}}$ and is actually an atom of $\mathcal{N}_{m}$ on $]s,t]\times \mathcal{P}^{0}_{m}\setminus \{0_{\cro{m}}\}$.
We then have the following identity:
$$\Pi^{m}(s,t)=\Pi^{n}_{|\cro{m}}(s,t).$$ 
This compatibility property allows us to define a unique process $(\Pi(s,t), -\infty<s\leq t<\infty)$ such that for all $s\leq t$, $\Pi_{|\cro{n}}(s,t)=\Pi^{n}(s,t)$. The collection $(\Pi(s,t), -\infty<s\leq t<\infty)$ is by construction a flow in the sense of Definition 1.
Obviously, the process $(\Pi(t), t\geq 0):=(\Pi(0,t), t\geq 0)$ is a standard distinguished coalescent with coagulation measure $\mu$. Interpreting a classical coagulation as a reproduction and a coagulation with the distinguished block as an immigration event, a coagulation measure $\mu$ may be viewed as encoding the births and the immigration rates in some population when time goes forward.
\\

In the same vein as Bertoin and Le Gall's flows, a population model is embedded in the dual flow, as we shall see. 
\subsection{The dual flow}
Let $\hat{\mathcal{N}}$ be the image of $\mathcal{N}$ by the time reversal $t\mapsto -t$ and consider the filtration $(\hat{\mathcal{F}}_{t},t\geq 0):=(\sigma (\hat{\mathcal{N}}_{[0,t]\times \mathcal{P}^{0}_{\infty}}),t\geq 0).$ For all $s\leq t$, we denote by $\hat{\Pi}(s,t)$ the partition $\Pi(-t,-s)$. The process $(\hat{\Pi}(s,t), -\infty <s\leq t <\infty)$ is called the dual flow. By a slight abuse of notation the process $(\hat{\Pi}(0,t), t\geq 0)$ will be denoted by $(\hat{\Pi}(t), t\geq 0)$. Plainly the following cocycle property holds for every $s\geq 0$, $$\hat{\Pi}(t+s)=coag(\hat{\Pi}(t,t+s), \hat{\Pi}(t)).$$
We stress that the partition $\hat{\Pi}(t,t+s)$ is exchangeable, independent of $\hat{\mathcal{F}}_{t}$ and has the same law as $\hat{\Pi}(s)$. The cocycle property yields immediately that $(\hat{\Pi}(t), t\geq 0)$ and its restrictions $(\hat{\Pi}_{|\cro{n}}(t), t\geq 0)$ are Markovian. The following property ensures that $(\hat{\Pi}(t), t\geq 0)$ is actually strongly Markovian.
\begin{prop}
The semigroup of the process $(\hat{\Pi}(t), t\geq 0)$ verifies the Feller property. For any continuous map $\phi$ from $\mathcal{P}^{0}_{\infty}$ to $\mathbb{R}$, the map $\pi \mapsto \mathbb{E}[\phi(coag(\hat{\Pi}(t), \pi))]$ is continuous and $\mathbb{E}[\phi(coag(\hat{\Pi}(t), \pi))]\underset{t \rightarrow 0}{\rightarrow}\phi(\pi).$
\end{prop}
\textit{Proof.} This is readily obtained thanks to the continuity of coagulation maps.$\square$
\begin{remark}
We stress that the process $(\hat{\Pi}(t), t \geq 0)$ is not a coalescent process. However, since the Poisson random measures $\mathcal{N}$ and $\hat{\mathcal{N}}$ have the same law, the process $(\hat{\Pi}(t), t \geq 0)$ has the same one-dimensional marginals as a standard coalescent with coagulation measure $\mu$.
\end{remark}
As explained in Section 2.1, the partition-valued process $(\hat{\Pi}(t), t\geq 0)$ may be viewed as a population model forward in time. For every $t\geq s\geq 0$ and every $k\in \mathbb{N}$, we shall interpret the block $\hat{\Pi}_{k}(s,t)$ as the descendants at time $t$ of the individual $k$ living at time $s$ and the distinguished block $\hat{\Pi}_{0}(s,t)$ as the descendants of the generic immigrant. Thanks to the cocycle property, for all $0 \leq s \leq t$, we have
$$\hat{\Pi}_{k}(t)=\bigcup_{j\in \hat{\Pi}_{k}(s)}\hat{\Pi}_{j}(s,t)$$
and the ancestor living at time $s$ of any individual $j$ at time $t$ is given by $\alpha_{\hat{\Pi}(s,t)}(j)$. The random distinguished partition $\hat{\Pi}(t)$ is exchangeable and possesses asymptotic frequencies. For all $i\geq 0$, we shall interpret $|\hat{\Pi}_{i}(t)|$ as the fraction of the population at time $t$ which is descendent from $i$. We stress that as in Donnelly-Kurtz's construction \cite{DonnKurtz} and the generalisation \cite{Birk2}, the model is such that the higher the individual is, the faster his descendants will die. Namely, for all $t\geq 0$ and all $j\geq 0$ we have $\alpha_{\hat{\Pi}(t)}(j)\leq j$ and then for all individuals $i<j$, the descendants of $i$ will always extinct after that of $j$.
\begin{remark}
When neither simultaneous multiple births nor immigration is taken into account, the measure $\mu$ is carried on the simple partitions (meaning with only one non-trivial block) with a distinguished block reduced to $\{0\}$ and we recover exactly the "modified" lookdown process of Donnelly-Kurtz for the $\Lambda$-Fleming-Viot process (see p195-196 of \cite{DonnKurtz}) . Moreover, the partition-valued process $(\hat{\Pi}(t), t\geq 0)$ corresponds in law with those induced by the dual flow of Bertoin and Le Gall $(\hat{B}_{t},t\geq 0)$ using the paint-box scheme, see \cite{LGB1}. 
\end{remark}
Let us study the genealogical process of this population model. We will recover a distinguished coalescent. Let $T>0$ be a fixed time and consider the population at time $T$. By definition, the individuals $k$ and $l$ have the same ancestor at time $T-t$ if and only if $k$ and $l$ belong to the same block of $\hat{\Pi}(T-t,T)$. Moreover by definition of the dual flow, $$(\hat{\Pi}(T-t,T), t\in [0,T])=(\Pi(-T,-T+t), t\in [0,T])$$ which is a distinguished coalescent with coagulation measure $\mu$ on the interval $[0,T]$.
\\
\\
In the same way as Donnelly and Kurtz in \cite{DonnKurtz}, we associate initially to each individual a \textit{type} represented by a point in a metric space $E\cup\{\partial\}$, where $\partial$ is an extra point not belonging to $E$ representing the distinguished type of the immigrants. The choice of $E$ does not matter in our setting and for the sake of simplicity, we choose for $E$ the interval $]0,1]$, and for distinguished type $\partial=0$. The generic external immigrant $0$ has the type $U_{0}:=0$. At any time $t$, each individual has the type of its ancestor at time $0$. In other words, for any $k\in \mathbb{N}$ the type of the individual $k$ at time $t$ is $U_{\alpha_{\hat{\Pi}(t)}(k)}$. The exchangeability properties of Section 2.2 will allow us to define and characterize the generalized Fleming-Viot process with immigration.
\begin{remark} We assume in this work only one source of immigration but several sources may be considered by distinguishing several blocks and types. In this work no mutation is assumed on the types. Birkner \textit{et al.} defined in \cite{Birk2} a generalization of the lookdown representation and the $\Xi$-Fleming-Viot process with mutations. Assuming that neither immigration nor mutation is taken into account, we will recover a process with the same law as a $\Xi$-Fleming-Viot process (the measure $\Xi$ is defined by $\Xi:=c_{1}\delta_{0}+\sum_{i\geq 1}s_{i}^{2}\nu(ds)$). However, we stress that our method differs from that of Birkner et al.  
\end{remark}
\section{Generalized Fleming-Viot processes as de Finetti measures}
We define and characterize in this section a measure-valued process which represents the frequencies of the types in the population at any time. This process will be called the generalized Fleming-Viot process with immigration and will be explicitly related to the \textit{forward partition} process $(\hat{\Pi}(t), t\geq 0)$. In the same way as in \cite{LGB1} and \cite{Birk2}, a duality argument allows us to characterize in law the GFVIs.
\subsection{Generalized Fleming-Viot processes with immigration} 
Let $\rho \in \mathcal{M}_{1}$, we assume that the initial types $(U_{i}, i\geq 1)$ are i.i.d. with law $\rho$ and independent of $\hat{\mathcal{N}}$. For all $t\geq 0$, the random partition $\hat{\Pi}(t)$ is exchangeable and applying Lemma \ref{exch}, we get that the sequence $(U_{\alpha_{\hat{\Pi}(t)}(l)}, l\geq 1)$ is exchangeable. We denote by $Z_{t}$ its de Finetti measure. Lemma \ref{exch} leads us to the following definition.
\begin{Def} \label{main}
The process $(Z_{t},t\geq 0)$ defined by $$Z_{t}:=|\hat{\Pi}_{0}(t)|\delta_{0}+\sum_{i=1}^{\infty}|\hat{\Pi}_{i}(t)|\delta_{U_{i}}+(1-\sum_{i=0}^{\infty}|\hat{\Pi}_{i}(t)|)\rho,$$ starting from $Z_{0}=\rho$, is called the generalized Fleming-Viot process with immigration.
\end{Def}
\begin{remark}For all $t\geq 0$, the random variable $Z_{t}$ can be viewed as the Stieltjes measure of a distinguished bridge (see \cite{coaldist}). Definition \ref{main} yields a paint-box representation of the population forward in time in the same vein as the dual flow of bridges of Bertoin and Le Gall.
\end{remark}
\begin{prop}
The process $(Z_{t}, t\geq 0)$ is Markovian with a Feller semigroup.
\end{prop}
\textit{Proof.} The sequence $(U_{\alpha_{\hat{\Pi}(s+t)}(l)},l\geq 1)$ is exchangeable with de Finetti measure $Z_{s+t}$. By the cocycle property of the dual flow, we have $\hat{\Pi}(s+t)=coag(\hat{\Pi}(s,s+t),\hat{\Pi}(s))$. Therefore, for all $l\geq 1$,
$$U_{\alpha_{\hat{\Pi}(s+t)}(l)}=U_{\alpha_{\hat{\Pi}(s)}\circ \alpha_{\hat{\Pi}(s,s+t)}(l)}.$$
By Lemma \ref{exch}, we immediately get that for all $t\geq 0$ and $s\geq 0$
$$Z_{s+t}=(1-\sum_{j\geq 0}|\hat{\Pi}_{j}(s,s+t)|)Z_{s}+\sum_{j\geq 1}|\hat{\Pi}_{j}(s,s+t)|\delta_{U_{\alpha_{\hat{\Pi}(s)}(j)}}+|\hat{\Pi}_{0}(s,s+t)|\delta_{0}.$$
We recall that $\hat{\Pi}(s,s+t)$ is independent of $\hat{\mathcal{F}}_{s}$ and hence of $Z_{s}$. By Theorem \ref{main}, conditionally on $Z_{s}$, the variables $(U_{\alpha_{\hat{\Pi}(s)}(j)} ,j\geq 1)$ are i.i.d. with distribution $Z_{s}$. The process $(Z_{t}, t\geq 0)$ is thus Markovian and its semigroup denoted by $R_{t}$ can be described as follows. For every $\rho \in \mathcal{M}_{1}$, $R_{t}(\rho,.)$ is the law of the random probability measure 
$$(1-\sum_{i\geq 0}|\hat{\Pi}_{i}(t)|)\rho+ \sum_{i\geq 1}|\hat{\Pi}_{i}(t)|\delta_{U_{i}}+|\hat{\Pi}_{0}(t)|\delta_{0}$$
where the variables $(U_{i}, i\geq 1)$ are i.i.d. distributed according to $\rho$ and independent of $\hat{\Pi}(t)$.
\\
We then verify that $R_{t}$ enjoys the Feller property. If $f$ is a continuous function from $\mathcal{M}_{1}$ to $\mathbb{R}$, the convergence in probability  when $t \rightarrow 0$ of $\hat{\Pi}(t)$ to $0_{\cro{\infty}}$ implies that $|\hat{\Pi}(t)|^{\downarrow}$ tends to $0$. We then have the convergence of $R_{t}f$ to $f$ when $t\rightarrow 0$. Plainly, for any sequence $(\rho^{n}, n\geq 1)$ which weakly converges to $\rho$, $R_{t}(\rho^{n},.)$ converges to $R_{t}(\rho,.)$. The Feller property is then established. $\square$
\subsection{Infinitesimal generator, core and martingale problem}
As in the articles of Bertoin-Le Gall \cite{LGB1}, Donnelly-Kurtz \cite{DonnKurtz} and Birkner et al \cite{Birk2}, the characterization in law of a GFVI will be obtained by a duality argument. Let $f$ be a continuous function on $[0,1]^{p}$. Define a function from $\mathcal{M}_{1}\times \mathcal{P}^{0}_{p}$ to $\mathbb{R}$ by
$$\Phi_{f} : (\rho,  \pi) \in \mathcal{M}_{1}\times \mathcal{P}^{0}_{p} \mapsto \int_{[0,1]^{p+1}} \delta_{0}(dx_{0})\rho(dx_{1})...\rho(dx_{p})f(x_{\alpha_{\pi}(1)},...,x_{\alpha_{\pi}(p)}).$$
Let $(\Pi(t), t\geq 0)$ be a distinguished coalescent, then we have the following lemma, where the notations $\mathbb{E}^{\rho}$ and $\mathbb{E}^{\pi}$ refer respectively to expectation when $Z_{0}=\rho$ and $\Pi_{|\cro{p}}(0)=\pi$.
\begin{lem} \label{dual}
$$\mathbb{E}^{\rho}[\Phi_{f}(Z_{t},\pi)]=\mathbb{E}^{\pi}[\Phi_{f}(\rho, \Pi_{|\cro{p}}(t))].$$
\end{lem}
\textit{Proof of Lemma \ref{dual}.} Let $(U_{i}, i\geq 1)$ be independent and identically distributed with law $\rho$ and $U_{0}=0$ , we have
\begin{align*}
\mathbb{E}^{\rho}[\Phi_{f}(Z_{t},\pi)]&=\mathbb{E}^{\rho}[\int \delta_{0}(dx_{0})Z_{t}(dx_{1})...Z_{t}(dx_{p})f(x_{\alpha_{\pi}(1)},...,x_{\alpha_{\pi}(p)})]\\
&=\mathbb{E}[f(U_{\alpha_{\hat{\Pi}(t)}(\alpha_{\pi}(1))},...,U_{\alpha_{\hat{\Pi}(t)}(\alpha_{\pi}(p))})]\\
&=\mathbb{E}[\int \delta_{0}(dy_{0})\rho(dy_{1})...\rho(dy_{p})f(y_{\alpha_{coag(\pi, \hat{\Pi}(t))}(1)},...,y_{\alpha_{coag(\pi, \hat{\Pi}(t))}(p)})]\\
&=\mathbb{E}^{\pi}[\Phi_{f}(\rho, \Pi_{|\cro{p}}(t))].
\end{align*}
The second equality holds because $Z_{t}$ is the de Finetti measure of $(U_{\alpha_{\hat{\Pi}(t)}(i)},i\geq 1)$. Observing that $\alpha_{\hat{\Pi}(t)}\circ \alpha_{\pi}=\alpha_{coag(\pi, \hat{\Pi}(t))}$, we get the third equality. Moreover, for $t$ fixed the random partition $\hat{\Pi}(t)$ has the same law as a standard distinguished coalescent at time $t$ which yields the last equality.$\square$ 
\\
\\
The Kolmogorov equations ensure that the generator $\mathcal{L}$ of the process $(Z_{t},t\geq 0)$ verifies for all continuous functions $f$ on $[0,1]^{p}$, $\pi\in \mathcal{P}^{0}_{p}$ and $\rho \in \mathcal{M}_{1}$, 
\begin{align} \label{gene}
\mathcal{L} \Phi_{f}(., \pi)(\rho)&=\mathcal{L}_{p}^{*} \Phi_{f}(\rho,.)(\pi)
\end{align}
where $\mathcal{L}_{p}^{*}$ is the generator of the continuous time Markov chain $(\Pi_{|\cro{p}}(t), t\geq 0)$. The process $(Z_{t},t\geq 0)$ is then characterized in law by a triplet $(c_{0},c_{1},\nu)$ according to the decomposition of the coagulation measure $\mu$ given in Subsection \ref{coaldist}. Let $G_{f}$ be the map defined by 
$$G_{f}(\rho)=\int_{[0,1]^{p}}f(x_{1},...,x_{p})\rho(dx_{1})...\rho(dx_{p})=\Phi_{f}(\rho, 0_{\cro{p}}).$$ 
A classical way to characterize the law of a Fleming-Viot process is to study a martingale problem, see for example Theorem 3 of \cite{LGB1} or Proposition 5.2 of \cite{coaldist}. The following theorem claims that a generalized Fleming-Viot process with immigration solves a well-posed martingale problem. Let $f$ be a continuous function $f$ on $[0,1]^{p}$. According to (\ref{gene}), the operator $\mathcal{L}$ is such that
$$\mathcal{L}G_{f}(\rho)=\sum_{\pi \in \mathcal{P}^{0}_{p}} q_{\pi} \int_{[0,1]^{p}}[f(x_{\alpha_{\pi}(1)},...,x_{\alpha_{\pi}(p)})-f(x_{1},...,x_{p})]\delta_{0}(dx_{0})\rho(dx_{1})...\rho(dx_{p}).$$
\begin{thm} \label{MPB} The law of the process $(Z_{t},t\geq 0)$ is characterized by the following martingale problem. For every integer $p\geq 1$ and every continuous function $f$ on $[0,1]^{p}$, the process
\begin{center} $G_{f}(Z_{t})-\int_{0}^{t}\mathcal{L}G_{f}(Z_{s})ds$ 
\end{center}
is a martingale in $(\hat{\mathcal{F}}_{t}, t\geq 0)$.
\end{thm}
\begin{remark} When $\mu$ is supported on the simple distinguished partitions, we recover the well-posed martingale problem of Section 5, Lemma 5.2 of \cite{coaldist}. We then identify the processes obtained from stochastic flows of bridges in \cite{LGB1} and \cite{coaldist} and from stochastic flows of partitions.
\end{remark}
\textit{Proof.}
Using (\ref{gene}), and applying Theorem 4.4.2 in \cite{EthierKurtz}, we get that there is at most one solution to the martingale problem. Dynkin's formula implies that the process in the statement is a martingale.
$\square$
\\
The following theorem yields an explicit formula for the generator of $(Z_{t}, t\geq 0)$. Its proof differs from that given in \cite{Birk2} for the $\Xi$-Fleming-Viot process. 
\begin{thm} \label{generator} The infinitesimal generator $\mathcal{L}$ of $(Z_{t}, t\geq 0)$ has the following properties:
\\
\begin{itemize}
\item[(i)]
For every integer $p\geq 1$ and every continuous function $f$ on $[0,1]^{p}$, we have
$$\mathcal{L}G_{f}=\mathcal{L}^{c_{0}}G_{f}+\mathcal{L}^{c_{1}}G_{f}+\mathcal{L}^{\nu}G_{f}$$
where
\begin{align*}
\mathcal{L}^{c_{0}}G_{f}(\rho)&:=c_{0}\sum_{1\leq i\leq p}\int_{[0,1]^{p}}[f(\textsl{x}^{0,i})-f(\textsl{x})]\rho^{\otimes p}(d\textsl{x})\\
\mathcal{L}^{c_{1}}G_{f}(\rho)&:=c_{1}\sum_{1\leq i<j\leq p}\int_{[0,1]^{p}}[f(\textsl{x}^{i,j})-f(\textsl{x})]\rho^{\otimes p}(d\textsl{x})\\
\mathcal{L}^{\nu}G_{f}(\rho)&:=\int_{\mathcal{P}_{\textbf{m}}}\lbrace\mathbb{E}[G_{f}(\bar{s}\rho+s_{0}\delta_{0}+\sum_{i\geq 1}s_{i}\delta_{U_{i}})]-G_{f}(\rho)\rbrace\nu(ds)
\end{align*}
where $\textsl{x}$ denotes the vector $(x_{1},...,x_{p})$ and
\begin{itemize}
\item the vector $\textsl{x}^{0,i}$ is defined by $\textsl{x}^{0,i}_{k}=x_{k}$, for all $k\neq i$ and $\textsl{x}^{0,i}_{i}=0$,
\item the vector $\textsl{x}^{i,j}$ is defined by $\textsl{x}^{i,j}_{k}=x_{k}$, for all $k \neq j$ and $\textsl{x}^{i,j}_{j}=x_{i},$
\item the sequence $(U_{i},i\geq 1)$ is i.i.d. with law $\rho$ and $\bar{s}$ is the dust of $s$ meaning that $\bar{s}:=1-\sum_{i\geq 0}s_{i}$.
\end{itemize}
\item[(ii)] Let $\mathcal{D}$ stand for the domain of $\mathcal{L}$. The vector space generated by functionals of the type $G_{f}$ forms a core of $(\mathcal{L}, \mathcal{D})$.\\
\end{itemize}
\end{thm}
\textit{Proof.} (i) We have $$\mathcal{L} G_{f}(\rho)=\mathcal{L}_{p}^{*}\Phi_{f}(\rho, .)(0_{\cro{p}}).$$ 
Therefore, $$\mathcal{L} G_{f}(\rho)=\sum_{\pi \in \mathcal{P}^{0}_{p}}q_{\pi}[\Phi_{f}(\rho,\pi)-\Phi_{f}(\rho,0_{\cro{p}})]$$ with $q_{\pi}=\mu(\mathcal{P}^{0}_{\infty, \pi})$. The decomposition of $\mu$ with the triplet $(c_{0},c_{1},\nu)$ implies that 
\\
\begin{align*}
\underset{\pi \in \mathcal{P}^{0}_{p}} {\sum}q_{\pi}[\Phi_{f}(\rho,\pi)-\Phi_{f}(\rho,0_{\cro{p}})]=&c_{0}\underset{1\leq i\leq p}\sum[f(\textsl{x}^{0,i})-f(\textsl{x})]\rho^{\otimes p}(d\textsl{x})+c_{1}\sum_{1\leq i<j\leq p}[f(\textsl{x}^{i,j})-f(\textsl{x})]\rho^{\otimes p}(d\textsl{x})\\
+&\sum_{\pi\in \mathcal{P}^{0}_{p}}\int_{\mathcal{P}_{\textbf{m}}}\rho_{s}(\mathcal{P}^{0}_{\infty, \pi})\nu(ds)[\Phi_{f}(\rho, \pi)-\Phi_{f}(\rho, 0_{\cro{p}})].
\end{align*}
It remains to establish the following equality
$$\mathcal{L}^{\nu}G_{f}(\rho)=\sum_{\pi\in \mathcal{P}^{0}_{p}}\int_{\mathcal{P}_{\textbf{m}}}\rho_{s}(\mathcal{P}^{0}_{\infty, \pi})\nu(ds)[\Phi_{f}(\rho, \pi)-\Phi_{f}(\rho, 0_{\cro{p}})].$$
Let $s\in \mathcal{P}_{\textbf{m}}$, as already mentioned, we denote by $\bar{s}$ its dust. Let $(U_{i}, i\geq 1)$ be i.i.d. random variables with distribution $\rho$. Denoting by $\Pi$ an independent $s$-distinguished paint-box, the variables $(U_{\alpha_{\Pi}(j)}, j\geq 1)$ are exchangeable with a de Finetti measure which has the same law as 
$$\bar{s}\rho+\sum_{i\geq 1}s_{i}\delta_{U_{i}}+s_{0}\delta_{0}.$$
Thus, we get
$$\mathbb{E}[G_{f}(\bar{s}\rho+\sum_{i\geq 1}s_{i}\delta_{U_{i}}+s_{0}\delta_{0})]=\mathbb{E}[f(U_{\alpha_{\Pi}(1)},...,U_{\alpha_{\Pi}(p)})].$$
Moreover,
$$\mathbb{E}[f(U_{\alpha_{\Pi}(1)},...,U_{\alpha_{\Pi}(p)})]-\mathbb{E}[f(U_{1},...,U_{p})]=\sum_{\pi\in \mathcal{P}^{0}_{p}}\mathbb{P}[\Pi_{|\cro{p}}=\pi]\left(\mathbb{E}[f(U_{\alpha_{\pi}(1)},...,U_{\alpha_{\pi}(p)})]-\mathbb{E}[f(U_{1},...,U_{p})]\right).$$
By definition, $\mathbb{P}[\Pi_{|\cro{p}}=\pi]=\rho_{s}(\mathcal{P}^{0}_{\infty, \pi})$ and we get by integrating on $\mathcal{P}_{\textbf{m}}$:
\begin{align*}
\sum_{\pi\in \mathcal{P}^{0}_{p}}\int_{\mathcal{P}_{\textbf{m}}}\rho_{s}&(\mathcal{P}^{0}_{\infty, \pi})\nu(ds)\left(\mathbb{E}[f(U_{\alpha_{\pi}(1)},...,U_{\alpha_{\pi}(p)})]-\mathbb{E}[f(U_{1},...,U_{p})]\right)\\
&=\int_{\mathcal{P}_{\textbf{m}}}\mathbb{E}[f(U_{\alpha_{\Pi}(1)},...,U_{\alpha_{\Pi}(p)})-f(U_{1},...,U_{p})]\nu(ds)\\
&=\int_{\mathcal{P}_{\textbf{m}}}\mathbb{E}[G_{f}(\bar{s}\rho+\sum_{i\geq 1}s_{i}\delta_{U_{i}}+s_{0}\delta_{0})-G_{f}(\rho)]\nu(ds).
\end{align*}
Therefore, the statement of \textit{(i)} is obtained.
\\
\\
(ii) The previous calculation yields that the map $\rho \mapsto \mathcal{L}G_{f}(\rho)$ is a linear combination of functionals of the type $\Phi_{f}(\rho, \pi)$. Besides, for all $\pi\in \mathcal{P}^{0}_{p}$, the map $\rho \mapsto \Phi_{f}(\rho, \pi)$ can be written as $G_{g}(\rho)$ with $g$ the continuous function on $[0,1]^{\#\pi-1}$ defined by
\begin{center} $g(x_{1},...,x_{\#\pi-1})=f(x_{\alpha_{\pi}(1)},...,x_{\alpha_{\pi}(p)})$ with $x_{0}=0$.
\end{center} 
Therefore, denoting by $D$ the vector space generated by the functionals of type $G_{f}$, the space $D$ is invariant under the action of the generator $\mathcal{L}$. Considering the maps of the form $f(x_{1},...,x_{p})=g(x_{1})...g(x_{p})$, we get that $D$ contains the linear combinations of functionals $\rho\mapsto \langle g,\rho \rangle^{p}$. By the Stone-Weierstrass theorem, these functionals are dense in the space of continuous functions on $\mathcal{M}_{1}$. Thanks to the Feller property of $(Z_{t},t\geq 0)$ and according to Proposition 19.9 of \cite{Kallenberg}, the space $D$ is a core. Thus, the explicit expression of $\mathcal{L}$ restricted to $D$, given in the statement, determines the infinitesimal generator $\mathcal{L}$ of $(Z_{t},t\geq 0)$. $\square$
\section{Extinction of the initial types}
Let $(Z_{t}, t\geq 0)$ be a GFVI characterized in law by the triplet $(c_{0},c_{1},\nu)$. The extinction of the initial types corresponds to the absorption of $(Z_{t}, t\geq 0)$ in $\delta_{0}$. It means for the forward partition-valued process $(\hat{\Pi}(t),t\geq 0)$ to be absorbed at the trivial partition $(\{\mathbb{Z}_{+}\},\emptyset,...)$.
We are interested in this section to determine under which conditions on $(c_{0},c_{1},\nu)$ the extinction occurs almost surely. We are only able to give a sufficient condition. However this condition is also necessary when the measure $\nu$ satisfies an additional assumption of regularity (as in Limic's article \cite{Limic}). We stress that as mentioned in Proposition 5.2 and Remark 5.3 of \cite{Birk2}, the absorption of a $\Xi$-Fleming-Viot process (without immigration) is closely related to the coming down from infinity of the $\Xi$-coalescent.
\subsection{Extinction criterion}
Define the following subspace of $\mathcal{P}_{\textbf{m}}$
\begin{center} $\mathcal{P}^{f}_{\textbf{m}}:=\{s\in \mathcal{P}_{\textbf{m}}; \sum_{i=0}^{n}s_{i}=1$ for some finite $n\}.$\end{center}
As in Section 5.5 of Schweinsberg's article \cite{Schweinsberg}, we consider the following cases:
\begin{itemize}
\item Suppose $\nu(\mathcal{P}^{f}_{\textbf{m}})=\infty$ then by a basic property of Poisson random measure, we know that $T_{f}:=\inf \{t>0; (t, \pi)$ is an atom of $\hat{\mathcal{N}}$ such that $\#\pi<\infty\}=0$ almost surely. We deduce that immediately after $0$, there is only a finite number of types and then the extinction occurs almost surely in a finite time.
\item Suppose $0<\nu(\mathcal{P}^{f}_{\textbf{m}})<\infty$ then $T_{f}$ is exponential with parameter $\nu(\mathcal{P}^{f}_{\textbf{m}})$. At time $T_{f}$ only a finite number of types reproduces and extinction will occur almost surely. This allows us to reduce the problem to the case when $\nu(\mathcal{P}^{f}_{\textbf{m}})=0$.
\end{itemize}
Suppose henceforth that $\nu(\mathcal{P}^{f}_{\textbf{m}})=0$. We define $$\zeta(q):=c_{1}q^{2}/2+ 
\int_{\mathcal{P}_{\textbf{m}}}\left(\sum_{i=1}^{\infty}(e^{-qs_{i}}-1+qs_{i})\right)\nu(ds).$$
\begin{thm} \label{ext}
If the following conditions are fulfilled:
\begin{itemize}
\item[i)] $c_{0}+\nu(s\in \mathcal{P}_{\textbf{m}}, s_{0}>0)>0$
\item[ii)] $\int_{a}^{\infty}\frac{dq}{\zeta(q)}<\infty$ for some $a>0$ (and then automatically for all $a>0$)
\end{itemize}
then the generalized Fleming-Viot process with immigration $(Z_{t},t\geq 0)$ is absorbed in $\delta_{0}$ almost surely.
\\
\\
Moreover, under the regularity condition $(R)$: $$\int_{\mathcal{P}_{m}}(\sum_{i=0}^{\infty}s_{i})^{2}\nu(ds)<\infty,$$
the conditions i) and ii) are necessary. 
\end{thm}
Following the terminology of Limic in \cite{Limic}, the assumption $(R)$ is called the \textit{regular} case. 
\begin{remark}
Some weaker assumption than $(R)$ may be found under which the conditions i) and ii) are still necessary (see Schweinsberg's article \cite{Schweinsberg} for the case of $\Xi$-coalescents). For sake of simplicity, we will not focus on that question in this article.
\end{remark} 
Theorem \ref{ext} extends Theorem 4.1 in \cite{coaldist}. Namely, the $M$-generalized Fleming-Viot processes always verify $(R)$. Moreover assume that $(R)$ holds and that the immigration and the reproduction never happen simultaneously. Therefore, the measure $\nu$ can be decomposed in the following way: $$\nu=\nu 1_{\{s\in \mathcal{P}_{\textbf{m}}; s=(s_{0},0,...)\}}+\nu 1_{\{s\in \mathcal{P}_{\textbf{m}}; s_{0}=0\}}.$$
Let us define
\begin{itemize}
\item $\nu_{0}:=\nu 1_{\{s\in \mathcal{P}_{\textbf{m}}; s=(s_{0},0,...)\}}$ which is viewed as a measure on $[0,1]$, encoding immigration rate, such that $\int_{0}^{1}s_{0}\nu_{0}(ds_{0})<\infty$.
\item $\nu_{1}:=\nu 1_{\{s\in \mathcal{P}_{\textbf{m}}; s_{0}=0\}}$ such that $\int_{\mathcal{P}_{\textbf{m}}}\left(\sum_{i\geq 1}s^{2}_{i}\right)\nu_{1}(ds)< \infty$
\end{itemize}
The map $\zeta$ does not depend on $\nu_{0}$. We deduce that, in this setting, the immigration has no impact on the extinction occurrence. However, we stress that the measure $\nu$ may be carried on $\{s\in \mathcal{P}_{\textbf{m}}; s_{0}>0, s_{1}>0\}$. 
\begin{rem} 
If the measure $\nu$ is carried on the set $\Delta:=\{(s_{0},s_{1})\in [0,1]^{2}; s_{0}+s_{1}<1\}$ then the regularity assumption $(R)$: $\int_{\Delta} (s_{0}+s_{1})^{2}\nu(ds)<\infty$ is still satisfied. Indeed for all $(s_{0},s_{1})\in \Delta$, $(s_{0}+s_{1})^{2}\leq 3s_{0}+ s_{1}^{2}$. In this setting, Theorem \ref{ext} gives a necessary and sufficient condition for extinction.\end{rem}
By duality, the extinction occurs if and only if there is an immigration and the embedded distinguished coalescent $(\Pi(t),t\geq 0)$ comes down from infinity, meaning that its number of blocks becomes instantaneously finite. We shall investigate this last question in the following subsection.
\subsection{Proof of Theorem \ref{ext}: coming down from infinity for a distinguished coalescent}
Let $(\Pi(t), t\geq 0)$ be a distinguished coalescent with triplet $(c_{0},c_{1},\nu)$.
\begin{Def} We say that a distinguished coalescent \textit{comes down from infinity} if \begin{center} $\mathbb{P}(\#\Pi(t)<\infty$ for all $t>0)=1$. \end{center}
\end{Def}
In the same manner as Schweinsberg's article \cite{Schweinsberg} Section 5.5, we will get a sufficient condition which will be necessary in a so-called \textit{regular} case (in the same sense as in Limic's article \cite{Limic}). The arguments used in \cite{coaldist} to study the coming down for an $M$-coalescent may be adapted in this more general framework. 
\begin{lem} \label{coming}
Consider the process $(\#\Pi(t),t\geq 0)$ of the number of blocks in a distinguished exchangeable coalescent. Provided that $\nu(\mathcal{P}^{f}_{\textbf{m}})=0$, there are two possibilities for the evolution of the number of blocks in $\Pi$: either $\mathbb{P}(\#\Pi(t)=\infty$ for all $t\geq 0)=1$ or $\mathbb{P}(\#\Pi(t)<\infty$ for all $t>0)=1$.
\end{lem}
\textit{Proof.} This is an easy adaptation of Lemma 31 in \cite{Schweinsberg}. $\square$
\\
\\
Let $\pi\in \mathcal{P}^{0}_{n}$, with $\#(\pi_{0}\setminus \{0\})=k_{0}$, $\#\pi_{1}=k_{1},..., \#\pi_{r}=k_{r}$ where $r\geq 0, k_{0}\geq 0$ and $k_{i}\geq 2$ for all $i\in [r]$ and $\sum_{i=0}^{r}k_{i}\leq n.$
We will denote by $\lambda_{n,k_{0},...,k_{r}}$ the jump rate $q_{\pi}$. The decomposition of the coagulation measure $\mu$ provides an explicit formula for $\lambda_{n,k_{0},...,k_{r}}$, however its expression is rather involved and we will not use it here. We stress that by exchangeability the quantity $\lambda_{n,k_{0},...,k_{r}}$ does not depend on the order of the integers $k_{1},...,k_{r}$.
From a partition with $n$ blocks, $k_{0}$-tuple, $k_{1}$-tuple,...,$k_{r}$-tuple merge simultaneously with rate $\lambda_{n,k_{0},...,k_{r}}$. The $k_{0}$-tuple represents the blocks coagulating with the distinguished one.
\\
\\
Define $N(n,k_{0},...,k_{r})$ to be the number of different simultaneous choices of a $k_{0}$-tuple, a $k_{1}$-tuple,.., and a $k_{r}$-tuple from a set of $n$ elements with $k_{0}\geq 0, k_{i}\geq 2$ for $i\in [r]$. The exact expression may be found but is not important in the rest of the current analysis. Denoting by $\Pi^{*}(t)=(\Pi_{1}(t),...,)$, we determine the generator $G^{[n]}$ of $(\#\Pi^{*}_{|[n]}(t),t\geq 0)$. 
Let $f$ be any map from $[n]$ to $\mathbb{R}$, $$G^{[n]}f(l)=\sum_{r=0}^{\lfloor l/2 \rfloor}\sum_{\substack{k_{0},\{k_{1},...k_{r}\};\\ \sum_{i=0}^{r}k_{i}\leq l}}N(l,k_{0},...,k_{r})\lambda_{l,k_{0},...,k_{r}}[f(l-(k_{0}+...+k_{r})+r)-f(l)]$$
As already mentioned, for each fixed $k_{0}$, we do not have a separate term for each possible ordering of $k_{1},...,k_{r}$. That is why the inner sum extends over $k_{0}\geq 0$ and the multiset $\{k_{1},...,k_{r}\}$ such that $\sum_{i=0}^{r}k_{i}\leq l$. 
\\
\\
We define the map $\Phi$ such that
$\Phi(n)$ is the total rate of decrease in the number of blocks in $(\Pi^{*}(t), t\geq 0)$ when the current configuration has $n$ blocks. We get
$$\Phi(n)=\sum_{r=0}^{\lfloor n/2 \rfloor}\sum_{\substack{k_{0},\{k_{1},...k_{r}\};\\ \sum_{i=0}^{r}k_{i}\leq n}}N(n,k_{0},...,k_{r})\lambda_{n,k_{0},...,k_{r}}[k_{0}+...+k_{r}-r].$$
\begin{lem} \label{concav}
\begin{itemize}
\item[i)] A more tractable expression for $\Phi$ is given by the following
$$\Phi(q)=c_{0}q+\frac{c_{1}}{2}q(q-1)+\int_{\mathcal{P}_{\textbf{m}}}(qs_{0}+\sum_{i=1}^{\infty}(qs_{i}-1+(1-s_{i})^{q})\nu(ds).$$
\item[ii)]  Define $$\Psi(q):=c_{0}q+c_{1}q^{2}/2+ \int_{\mathcal{P}_{\textbf{m}}}\left(qs_{0}+\sum_{i=1}^{\infty}(e^{-qs_{i}}-1+qs_{i})\right)\nu(ds).$$
There exist $C$ and $C'$ two nonnegative constants such that $C\Psi(q)\leq \Phi(q)\leq C'\Psi(q)$. 
\item[iii)] The map $q\mapsto \Psi(q)/q$ is concave.
\end{itemize}
\end{lem}
\textit{Proof.} Proof of i): We have $N(q,1)=\binom{q}{1}=q$ and $N(q,0,2)=\binom{q}{2}$. Using the binomial formula, the first two terms are plain. We focus now on the integral term. Let $s\in \mathcal{P}_{\textbf{m}}$ and $\pi$ be a $s$-distinguished paint-box (see Definition 2.3 in \cite{coaldist}). To compute the total rate of decrease in the number of blocks from a configuration with $q\geq 1$ blocks, let us consider $$Y^{(q)}_{l}(\pi):=\#\{k\in [q]; \alpha_{\pi}(k)=l\}.$$ 
Conditionally given $|\pi_{l}|$, the variable $Y^{(q)}_{l}(\pi)$ is binomial with parameters $(|\pi_{l}|,q)$ (degenerated in the case of $|\pi_{l}|=0$) and we have 
\begin{align*}
\Phi(q)&=c_{0}q+\frac{c_{1}}{2}q(q-1)+\int_{\mathcal{P}_{m}}\left(\mathbb{E}[Y^{(q)}_{0}(\pi)]+\sum_{l=1}^{\infty}\mathbb{E}[Y^{(q)}_{l}(\pi)-1_{\{Y^{(q)}_{l}(\pi)>0\}}]\right)\nu(ds)\\
&=c_{0}q+\frac{c_{1}}{2}q(q-1)+\int_{\mathcal{P}_{m}}\left(qs_{0}+\sum_{l=1}^{\infty}[qs_{l}-1+(1-s_{l})^{q}]\right)\nu(ds).
\end{align*} 
Note that these computations are exactly the same as those pages 224-225 in \cite{Limic} using the "coloring procedure".
\\
Proof of ii) Same calculations as in Lemma 8 in Limic's article \cite{Limic}.\\
Proof of iii) We remark that $\Psi$ is the Laplace exponent of a spectrally positive L\'evy process. Therefore the map $h:q\mapsto \Psi(q)/q$ is the Laplace exponent of a subordinator which is concave. 
$\square$
\\
\\
The following theorem may be compared with Schweinsberg's criterion in \cite{CDI} and Theorem 4.1 in \cite{coaldist}.
\begin{thm}
The convergence of the series $\sum_{n\geq 1}\frac{1}{\Phi(n)}$ implies the coming down from infinity of the distinguished coalescent.
\\
\\
Under the regularity condition $(R)$: $$\int_{\mathcal{P}_{m}}(\sum_{i=0}^{\infty}s_{i})^{2}\nu(ds)<\infty,$$
the convergence of the series is necessary. 
\end{thm}
We refer to Example 34 p40 in \cite{Schweinsberg} and Section 3.2 p231 in \cite{Limic} for a coalescent which is not regular, comes down from infinity, with a divergent series. We shall follow the proof of Theorem 4.1 in \cite{coaldist} and study some (super)martingales.
\begin{lem} Let $(\Pi(t),t\geq 0)$ be a distinguished coalescent with triplet $(c_{0},c_{1},\nu)$. Assume that $\nu(\mathcal{P}^{f}_{\textbf{m}})=0$. Let us define the fixation time $$\zeta:=\inf\{t\geq 0: \Pi(t)=\{\mathbb{Z}_{+},\emptyset,...\}\}.$$ 
The expectation of fixation time is bounded by $$\mathbb{E}[\zeta]\leq \sum_{n=1}^{\infty}1/\Phi(n).$$
As a consequence, if the series in the right-hand side converges, the fixation time is finite with probability one.
\end{lem}
\textit{Proof.} Assuming the convergence of the sum $\sum_{n=1}^{\infty}1/\Phi(n)$, we define $$f(l)=\sum_{k=l+1}^{\infty}1/\Phi(k).$$ It is easy to check directly from i) in Lemma \ref{concav} that the map $\Phi$ is increasing (see alternatively Lemma 28 in \cite{Schweinsberg}). We thus have 
$$f(l-(k_{0}+...+k_{r})+r)-f(l)=\sum_{k=l-(k_{0}+...+k_{r})+r+1}^{l}\frac{1}{\Phi(k)}\geq \frac{k_{0}+...+k_{r}-r}{\Phi(l)}.$$ 
Therefore $$G^{[n]}f(l)\geq \sum_{r=0}^{\lfloor l/2 \rfloor}\sum_{\substack{k_{0},\{k_{1},...,k_{r}\}\\ \sum_{i=0}^{r}k_{i}\leq l}}\lambda_{l,k_{0},...,k_{r}}N(l,k_{0},...,k_{r})\frac{k_{0}+...+k_{r}-r}{\Phi(l)}=1.$$ 
The process $f(\#\Pi^{*}_{|\cro{n}}(t))-\int_{0}^{t}G^{[n]}f(\#\Pi^{*}_{|\cro{n}}(s))ds$ is a martingale. The quantity $$\zeta_{n}:=\inf\{t\geq 0: \#\Pi^{*}_{|\cro{n}}(t)=0\}$$ is a finite stopping time. Let $k\geq 1$, applying the optional sampling theorem to the bounded stopping time $\zeta_{n}\wedge k$, we get : $$\mathbb{E}[f(\#\Pi^{*}_{|\cro{n}}(\zeta_{n}\wedge k))]-\mathbb{E}[\int_{0}^{\zeta_{n}\wedge k}G^{[n]} f(\#\Pi^{*}_{|\cro{n}}(s))ds]=f(n)$$
With the inequality $G^{[n]}f(l)\geq 1$, we deduce that $$\mathbb{E}[\zeta_{n}\wedge k]\leq \mathbb{E}[f(\#\Pi^{*}_{|\cro{n}}(\zeta_{n}\wedge k))]-f(n).$$ By monotone convergence and Lebesgue's theorem, we have $\mathbb{E}[\zeta_{n}]\leq f(0)-f(n)$.  Passing to the limit in $n$, we have $\zeta_{n} \uparrow \zeta_{\infty}:= \inf \{t; \#\Pi(t)=1\}$ and $f(n) \longrightarrow 0$, thus $$\mathbb{E}[\zeta_{\infty}]\leq f(0)=\sum_{k=1}^{\infty}1/\Phi(k).$$ $\square$
\\
\\
To prove that convergence of the series is necessary for coming down from infinity we follow the same steps as in Section 6 of \cite{coaldist}. Assuming that $(R)$ holds, that the coalescent comes down from infinity and that the series is infinite, we may define a supermartingale (thanks to Lemmas \ref{binom}, \ref{half} and \ref{superm} below) and find a contradiction by applying the optional sampling theorem (Lemma \ref{CN}). The proofs of these lemmas are easy adaptations of those in Section 6 of \cite{coaldist}. We simply give their statements and the corresponding references in \cite{coaldist}.
\\
\\
The following technical lemma allows us to estimate the probability for the sum of $n$ independent binomial variables to be larger than $n/2$. 
\begin{lem} \label{binom}
Let $s\in \mathcal{P}_{\textbf{m}}$. Let $\pi$ be an $s$-distinguished paint-box and the variables $Y^{(n)}_{l}(\pi)$ defined as in Lemma \ref{concav}. For every $n_{0}\geq 4$, provided that $\sum_{i=0}^{\infty}s_{i}$ is sufficiently small, there is the bound
$$\mathbb{P}[\exists n\geq n_{0}; \sum_{l=0}^{n}Y^{(n)}_{l}(\pi)>\frac{n}{2}]\leq \frac{\exp{(-n_{0}f(s))}}{1-\exp{(-f(s))}}$$
with $$f(s)=\frac{1}{2}\log(\frac{1}{\sum_{i\geq 0}s_{i}})-\sum_{l=0}^{\infty}\log\left(\frac{s_{l}}{\sum_{i\geq 0}s_{i}}+1-s_{l}\right).$$
\end{lem}
\textit{Proof}. Easy adaptation of the arguments of Lemma 6.2 in \cite{coaldist}.
$\square$
\begin{lem} \label{half} Assume that the coalescent comes down from infinity. With probability one, we have
$$\tau:=\inf\{t>0, \#\Pi(t)<\frac{\#\Pi(t-)}{2}\}>0.$$ Moreover, if we define $\tau_{n}:=\inf\{t>0, \#\Pi_{|\cro{n}}(t)\leq \frac{\#\Pi_{|\cro{n}}(t-)}{2}\}$, then the sequence of stopping times $\tau_{n}$ converges to $\tau$ almost surely.
\end{lem}
\textit{Proof.} Using the assumption $(R)$ and the above Lemma \ref{binom}, we get $\mathbb{E}[\mathcal{N}(\{(t,\pi); t\leq 1; \exists n\geq n_{0}, \sum_{l=0}^{n}Y^{(n)}_{l}(\pi)>\frac{n}{2}\})]<\infty$. The same arguments as in Lemma 6.3 of \cite{coaldist} yield the statement. 
$\square$
\begin{lem} \label{superm}
Assuming $(R)$, the coming down from infinity and that $\sum_{n\geq 1} \frac{1}{\Phi(n)} =\infty$. There exists a constant $C>0$ such that for all $n\geq 1$, $(e^{-Ct}f(\#\Pi^{*}_{|[n]}(t)))_{t\leq \tau_{n}}$ is a non-negative supermartingale.
\end{lem}
\textit{Proof.} With the part iii) of Lemma \ref{concav} and the assumption $(R)$, this is an easy adaptation of the arguments of Lemma 6.4 in \cite{coaldist}.
$\square$
\begin{lem} \label{CN}
Under $(R)$, if $\sum_{n\geq 1} \frac{1}{\Phi(n)}=\infty$ then $\Pi$ does not come down from infinity.
\end{lem}
\textit{Proof.} Assume that the coalescent comes down from infinity. Exactly as in Lemma 6.5 in \cite{coaldist}, we may apply the optional sampling theorem to the previous supermartingale and find a contradiction.
$\square$
\\
\\
Thanks to part (ii) of Lemma \ref{concav}, we have the following equivalence
\[\sum_{n\geq 1}\frac{1}{\Phi(n)}<\infty \Longleftrightarrow \int_{a}^{\infty}\frac{dq}{\Psi(q)}<\infty.\]
In order to establish Theorem \ref{ext}, it suffices to show that $$\int_{a}^{\infty}\frac{dq}{\Psi(q)}<\infty \Longrightarrow \int_{a}^{\infty}\frac{dq}{\zeta(q)}<\infty.$$
Plainly, the quantity $\frac{q}{\zeta(q)}$ is bounded. Therefore for some constants $c$ and $C$, $$\Psi(q)=\zeta(q)(1+c\frac{q}{\zeta(q)})\leq C \zeta(q).$$
Finally, we get that the conditions i) and ii) of Theorem \ref{ext} are necessary under the assumption $(R)$.
\section*{Appendix}
We restate and prove here Lemma \ref{exch}.\\
\\
\textbf{Lemma \ref{exch}} \textit{ Let $(U_{i},i\geq 1)$ be an infinite exchangeable sequence taking values in $[0,1]$, with de Finetti measure $\rho$ and fix $U_{0}=0$. Let $\pi$ be an independent exchangeable distinguished partition. Then the infinite sequence $(U_{\alpha_{\pi}(k)},k\geq 1)$ is exchangeable. Furthermore, its de Finetti measure is $$(1-\sum_{i\geq 0}|\pi_{i}|)\rho+ \sum_{i\geq 1}|\pi_{i}|\delta_{U_{i}}+|\pi_{0}|\delta_{0}.$$ }
\textit{Proof.} By de Finetti's theorem, without loss of generality we may directly assume that the sequence $(U_{i}, i\geq 1)$ is i.i.d. with a distribution $\rho \in \mathcal{M}_{1}$. We show that for all $n\geq 1$, the random vector $(U_{\alpha_{\pi}(1)},...,U_{\alpha_{\pi}(n)})$ is then exchangeable. Let $f_{1},...,f_{n}$ be $n$ measurable functions on $[0,1]$ and $[n]$ be the set $\{1,...,n\}$,
$$\mathbb{E}[f_{1}(U_{\alpha_{\pi}(1)})...f_{n}(U_{\alpha_{\pi}(n)})|\pi]=\left(\underset{i \in \pi_{0}\cap[n]}{\prod}f_{i}(0) \right) \left(\prod_{k\geq 1}\int_{0}^{1}\underset{i \in \pi_{k}\cap [n]}{\prod}f_{i}(u)\rho(du)\right).$$
Let $\sigma$ be a permutation of $\mathbb{Z}_{+}$ such that $\sigma(0)=0$, and $\eta$ be a permutation such that $\sigma^{-1}(\pi_{i})=\sigma \pi_{\eta(i)}$. We stress that $\eta(0)=0$, and we have
\begin{align*}
\mathbb{E}[f_{1}(U_{\alpha_{\pi}(\sigma(1))})...f_{n}(U_{\alpha_{\pi}(\sigma(n))})|\pi]&=\left(\underset{i \in \sigma^{-1}(\pi_{0})\cap [n]}\prod f_{i}(0)\right) \left(\prod_{k\geq 1}\int_{0}^{1}\underset{i \in \sigma^{-1}(\pi_{k})\cap[n]} {\prod}f_{i}(u)\rho(du)\right)\\
&=\left(\underset{i \in \sigma \pi_{0}\cap[n]}\prod f_{i}(0) \right) \left( \prod_{k\geq 1}\int_{0}^{1}\underset{i \in \sigma\pi_{\eta(k)}\cap [n]}{\prod}f_{i}(u)\rho(du)\right)\\
&=\left(\underset{i \in \sigma \pi_{0}\cap[n]}{\prod}f_{i}(0)\right) \left(\prod_{k\geq 1}\int_{0}^{1}\underset{i \in \sigma \pi_{k}\cap[n]} {\prod}f_{i}(u)\rho(du)\right).
\end{align*}
Therefore
$$\mathbb{E}[f_{1}(U_{\alpha_{\pi}(\sigma(1))})...f_{n}(U_{\alpha_{\pi}(\sigma(n))})|\pi]=\mathbb{E}[f_{1}(U_{\alpha_{\sigma \pi}(1)})...f_{n}(U_{\alpha_{\sigma \pi}(n)})|\pi].$$
The exchangeability of the partition $\pi$ ensures that $$\mathbb{E}[f_{1}(U_{\alpha_{\sigma \pi}(1)})...f_{n}(U_{\alpha_{\sigma \pi}(n)})]=\mathbb{E}[f_{1}(U_{\alpha_{\pi}(1)})...f_{n}(U_{\alpha_{\pi}(n)})],$$
which allows us to conclude that $(U_{\alpha_{\pi}(1)},...,U_{\alpha_{\pi}(n)})$ is exchangeable.\\

We have now to identify the de Finetti measure of this exchangeable sequence. This is an easy adaptation of the proof of Lemma 4.6 of \cite{coursbertoin}. For the sake of completeness, we give the proof in detail. The exchangeability is given by Lemma \ref{exch}, and so to prove the statement it suffices, by de Finetti's theorem, to study the limit when $m\rightarrow \infty$ of 
$$\frac{1}{m}\sum_{j=1}^{m}\delta_{U_{\alpha_{\pi}(j)}}.$$
We denote by $S$ the set of singletons in the partition $\pi$. That is
$$S:=\bigcup_{i\geq 0, |\pi_{i}|=0}\pi_{i}.$$
In the sequel, we will denote by $S^{c}$ the complement of $S$. The paint-box structure of an exchangeable partition tells us that for all $t\geq 0$, $S$ is empty or infinite. We decompose the sum in the following way
\begin{align*}
\frac{1}{m}\sum_{j=1}^{m}\delta_{U_{\alpha_{\pi}(j)}}&=\frac{1}{m}\sum_{j\in S, j\in \cro{m}}\delta_{U_{\alpha_{\pi}(j)}}+\frac{1}{m}\sum_{j\in S^{c}, j\in \cro{m}}\delta_{U_{\alpha_{\pi}(j)}}\\
&=\#(S\cap \cro{m})/m \frac{1}{\#(S\cap \cro{m})}\sum_{j \in S; j\in \cro{m}}\delta_{U_{\alpha_{\pi}(j)}}+\frac{1}{m}\sum_{j \in S^{c}, j\in \cro{m}}\delta_{U_{\alpha_{\pi}(j)}}.
\end{align*}
By de Finetti's theorem, we deduce that the first term converges to $|S|\rho$. Let us decompose the second term
$$\frac{1}{m}\sum_{j\in S^{c}, j\in \cro{m}}\delta_{U_{\alpha_{\pi}(j)}}=\frac{1}{m}\sum_{i \in \cro{m}; \pi_{i}\subset S^{c}}\#(\pi_{i}\cap \cro{m})\delta_{U_{i}}.$$
From Fatou's lemma, we get for any measurable bounded function $f$ on $[0,1]$
$$\underset{m \rightarrow \infty}{\liminf} \sum_{i \in \cro{m}; \pi_{i}\subset S^{c}}\frac{1}{m}\#(\pi_{i}\cap \cro{m})f(U_{i}) \geq |\pi_{0}|f(0)+\sum_{i \geq 1; \pi_{i}\subset S^{c}} |\pi_{i}|f(U_{i}).$$
Therefore, 
$$\underset{m \rightarrow \infty}{\liminf} \langle f; \frac{1}{m} \sum_{j=1}^{m}\delta_{U_{\alpha_{\pi}(j)}} \rangle\geq \langle f; (1-\sum_{i\geq 0}|\pi_{i}|)\rho+|\pi_{0}|\ \delta_{0}+\sum_{i \geq 1} |\pi_{i}|\ \delta_{U_{i}}\rangle.$$
The last sum extends over all $i\geq 1$ because if $\pi_{i}$ is included in $S$, then the quantity $|\pi_{i}|$ is $0$. 
\\
\\
In order to study the $\limsup$, define for all $\eta>0$,
$$J(\eta):=\{j \in \mathbb{Z}_{+}; |\pi_{j}|\geq \eta \}.$$
This set is finite, and we have
$$\frac{1}{m} \underset{j \in S^{c}, j \in \cro{m}}\sum \delta_{U_{\alpha_{\pi}(j)}}=\frac{1}{m} \underset{j \in J(\eta), j \in \cro{m}}\sum \delta_{U_{\alpha_{\pi}(j)}}+\frac{1}{m} \underset{j \in S^{c}\setminus J(\eta), j \in \cro{m}}\sum \delta_{U_{\alpha_{\pi}(j)}}$$
The first sum extends over a finite set, so we can interchange the sum and the limit. For all $\eta>0$,
$$\frac{1}{m} \underset{j \in J(\eta), j \in \cro{m}}\sum f(U_{\alpha_{\pi}(j)}) \underset{ m\rightarrow \infty}\rightarrow \sum_{j; \ \pi_{j}\subset J(\eta)}|\pi_{j}|f(U_{j}).$$
Let us study the second sum, denoting by $C$ a constant such that $|f|\leq C$. We have for large $m$,
$$\frac{1}{m} \underset{j \in S^{c}\setminus J(\eta), j \in \cro{m}}\sum f(U_{\alpha_{\pi}(j)})\leq C\frac{1}{m}\# ( \underset{j \in S^{c}\setminus J(\eta)} \bigcup \pi_{j}\cap \cro{m}).$$
When $m\rightarrow \infty$, the boundary term converges to $|\bigcup_{j \in \ S^{c}\setminus J(\eta)} \pi_{j}|$.
We then have
$$\underset{m \rightarrow \infty}\limsup \frac{1}{m} \underset{j \in \ S^{c}\cap\cro{m}}\sum f(U_{\alpha_{\pi}(j)})\leq \sum_{j; \ \pi_{j}\subset J(\eta)}|\pi_{j}|f(U_{j})+C|\bigcup_{j \in S^{c}\setminus J(\eta)} \pi_{j}|.$$
By definition, we have $1-|S|-\sum_{j \in J(\eta)}|\pi_{j}| \underset{\eta \rightarrow 0} \rightarrow 0$ and therefore
$$|\bigcup_{j \in S^{c}\setminus J(\eta)} \pi_{j}|\underset{\eta \rightarrow 0} \rightarrow 0.$$
Obviously, $$\sum_{j; \ \pi_{j}\subset J(\eta)}|\pi_{j}|f(U_{j})\underset{\eta \rightarrow 0}\rightarrow \sum_{j\geq 0}|\pi_{j}|f(U_{j}).$$ Combining all these results, we get
$$\underset{m \rightarrow \infty}{\limsup} \langle f; \frac{1}{m} \sum_{j=1}^{m}\delta_{U_{\alpha_{\pi}(j)}} \rangle\leq \langle f; (1-\sum_{i\geq 0}|\pi_{i}|)\rho+|\pi_{0}|\ \delta_{0}+\sum_{i \geq 1} |\pi_{i}|\ \delta_{U_{i}}\rangle.$$
We then obtain the statement of the proposition.$\square$
\\
\\
\textbf{Acknowledgments}. I would like to thank my advisor Jean Bertoin for his very useful advice, and the anonymous referees for their careful reading and helpful suggestions.

\end{document}